\documentclass[leqno, myheadings, twoside]{amsart}
\usepackage{amsfonts}
\usepackage{amsmath, amsthm, amssymb, amscd, amsxtra,graphicx}
\usepackage{latexsym, amsfonts}
\usepackage{url}
\usepackage{texdraw}
\usepackage{epsfig}
\makeatletter \@addtoreset{equation}{section}

\makeatletter \renewcommand{\@biblabel}[1]{#1.}

\theoremstyle{remark}

\begin{document}
\title [Heat invariants of the perturbed polyharmonic Steklov problem] {Heat invariants of the perturbed polyharmonic Steklov problem}
\author{Genqian Liu}


\maketitle Department of Mathematics, Beijing Institute of
Technology,
 Beijing 100081, the People's Republic of China.
 \ \
E-mail address:  liugqz@bit.edu.cn

\vskip 0.46 true cm

\vskip 0.45 true cm

\vskip 15 true cm

\begin{abstract}    For a given bounded domain $\Omega$ with smooth boundary in a smooth Riemannian manifold $(\mathcal{M},g)$, we  establish a procedure to get all the coefficients of the asymptotic expansion of the trace of the heat kernel associated with the perturbed polyharmonic Dirichlet-to-Neumann operator $\Lambda_m$ ($m\ge 1$) as $t\to 0^+$. We also explicitly calculate the first four coefficients of this asymptotic expansion. These coefficients (i.e., heat invariants) provide precise information for the area and curvatures of the boundary $\partial \Omega$ in terms of the spectrum of the perturbed polyharmonic Steklov problem. In particular, when $m=1$ and $q\equiv 0$ our work  recovers  the previous corresponding results in \cite{PS} and \cite{Liu3}.
   \end{abstract}

\vskip 1.39 true cm

\section{ Introduction}

\vskip 0.45 true cm

Let $\mathcal{M}$ be an $(n+1)$-dimensional, smooth, complete Riemannian
manifold  with  metric tensor $g=(g_{jk})$, and let $\Omega$ be
 a bounded domain in $\mathcal{M}$ with smooth boundary $\partial \Omega$.
  Let $(-\Delta_g)^m, \, m=1,2,\cdots$, be a polyharmonic operator, where the Laplace-Beltrami operator associated with the metric $g$ is given in local coordinates by
 \begin{eqnarray} \label{1-1} \Delta_g u=\frac{1}{\sqrt{|g|}}\sum_{j,k=1}^{n+1}
   \frac{\partial}{\partial x_j} \left( \sqrt{|g|}\,g^{jk}
   \frac{\partial u}{\partial x_k}\right),\end{eqnarray}
 and $(g^{jk})$ is the inverse of the metric tensor
 $(g_{jk})$ and $|g| =\mbox{det}\,g$.
  Let $q\in L^{\infty} (\Omega)$ be a real valued potential. Here and in what
follows $H^s(\Omega)$ and $H^s (\partial \Omega), \, s\in {\Bbb R}^1$, are the standard $L^2$-based Sobolev spaces in $\Omega$ and its boundary $\partial \Omega$, respectively, and $\nu$ is the inward unit normal to the
boundary.
 For $u\in H^m(\Omega)$, we denote \begin{eqnarray*} \gamma u=(u\big|_{\partial \Omega}, ((-\Delta_g) u)\big|_{\partial \Omega},\cdots, ((-\Delta_g)^{m-1} u)\big|_{\partial \Omega}).\end{eqnarray*}
Throughout the paper, we will assume that $n\ge 2m$ and $0$ is not the spectrum of the operator
\begin{eqnarray*}  (-\Delta_g)^m +q: \, H_0^m (\Omega) \to H^{-m} (\Omega) =(H^m_0(\Omega))'.\end{eqnarray*}
 It follows from \cite{KU} that for any $\phi=(\phi_0, \cdots, \phi_{m-1})\in \prod_{j=0}^{m-1} H^{m-2j-1/2} (\partial \Omega)$, the Dirichlet problem,
\begin{eqnarray} \label{1.1} \left\{ \begin{array}{ll}  \big((-\Delta_g)^m +q\big)u =0 \;\; &\mbox{in}\;\; \Omega,\\
  \gamma u=\phi \;\; &\mbox{on}\;\; \partial \Omega,  \end{array}\right.\end{eqnarray}
has a unique solution  $u\in H^{m}(\Omega)$.
Introducing the Neumann-type trace operator $\tilde{\gamma}$ by
\begin{eqnarray*}&&\tilde{\gamma}\, :\,H^{m}(\Omega)\to \prod_{j=0}^{m-1} H^{m-(2j+1)-1/2} (\partial \Omega),\\
&& \tilde{\gamma} u= (\partial_\nu u\big|_{\partial \Omega}, \partial_\nu((-\Delta_g) u)\big|_{\partial \Omega},  \cdots,
  \partial_\nu ((-\Delta_g)^{m-1}u)\big|_{\partial \Omega}),\end{eqnarray*}
we define the perturbed polyharmonic Dirichlet-to-Neumann operator  $\Lambda_m$ by
  \begin{eqnarray*} \Lambda_m \phi= \tilde{\gamma} u,  \end{eqnarray*}
where  $u\in H^{m} (\Omega)$ is the solution to the Dirichlet problem (\ref{1.1}). It follows from \cite{KU} (see also, \cite{Gr2} and \cite{McL}) that $\Lambda_m$ is well defined.
Let us also introduce the graph $\mathcal {C}$ of the perturbed polyharmonic Dirichlet-to-Neumann operator $\Lambda_m$  as follows:
\begin{eqnarray} \label{1.3*} \mathcal{C} =\{(\gamma u, \tilde {\gamma} u )\big| u\in H^{m} (\Omega),\, \big((-\Delta_g u)^m +q\big)u=0 \;\; \mbox{in}\;\; \Omega\}.\end{eqnarray}
 Notice that when $m=1$ and $q\equiv 0$, we recover the standard
Dirichlet-to-Neumann map for the Schr\"{o}dinger equation, given by
\begin{eqnarray} \label{1.2}  H^{1/2} (\partial \Omega) \ni \phi \mapsto \partial_\nu u\big|_{\partial \Omega} \in H^{-1/2} (\partial \Omega),
\end{eqnarray}
where $\Delta_g u=0 \;\;\mbox{in}\;\; \Omega$ and $u=\phi$ on $\partial \Omega$.

The perturbed polyharmonic Dirichlet-to-Neumann operator is a self-adjoint,  elliptic
pseudodifferential operator.
  This  operator has been discussed as inverse boundary problems in \cite{KU}, \cite{FKSU}, \cite{KLU}, \cite{Gr2} and \cite{Isa}.
 For example, when $m=1$, the harmonic Dirichlet-to-Neumann operator is  related to the
 Calder\'{o}n problem  of determining the anisotropic conductivity of a body from current and
voltage measurements at its boundary (\cite{Cal}, \cite{LU}, \cite{Cha},  \cite{GLU}, \cite{IMS}\cite{SU1}, \cite{SU2}, \cite{AP}, \cite{Ale},  \cite{IMS}). Let us remark that
the areas of physics and geometry where such operators occur, include the study
of fluid mechanics, vibration problems, the Kirchhoff plate equation in the theory of elasticity, and the study of the
higher-order Paneitz-Branson operator in conformal geometry (see \cite{GGS} or \cite{KU}).
 Because $\partial \Omega$ is
compact, the spectrum of $\Lambda_m$ is nonpositive, discrete and
unbounded (see p.$\,$95 of \cite{Ban}, \cite{FSc}). The spectrum $\{\lambda_k\}_{k=1}^\infty$ of this operator
is just the perturbed polyharmonic Steklov spectrum of the domain $\Omega$.
More precisely, \begin{eqnarray*} \left\{ \begin{array}{ll}   \big((-\Delta_g)^m  + q\big)u_k
=0 \quad  \; & \mbox{in}\;\; \Omega, \\
 \tilde{\gamma} u_k=-\lambda_k (\gamma u_k), \; \quad &\mbox{on}\;\; \partial \Omega,\end{array} \right.\end{eqnarray*}
where $\gamma u_k$ is the (vector-valued) eigenfunction corresponding to the $k$-th perturbed polyharmonic Steklov eigenvalue $\lambda_k$. All the eigenfunctions $\{\gamma u_k\}$ form an orthogonal basis in $(L^2(\partial \Omega))^m$.
 The study of
the spectrum of $\Lambda_m$ with $m=1$ and  $q\equiv 0$ was initiated by Steklov in 1902 (see \cite{St}), and has been investigated extensively in the past a century in the literature (see, \cite{Pay}, \cite{KK}, \cite{FK}, \cite{Sa}, \cite{Ta2}, \cite{SU2}, \cite{Liu2}).
 We consider  the  following equation:
  \begin{eqnarray}  \label {1-5}\left\{ \begin{array}{ll} \frac{\partial u(t, x)}{\partial t}= {\Lambda}_m u(t,x)  \;\;
\quad \; &\mbox{in}\;\;  [0, +\infty)\times \partial\Omega, \\
 \gamma u(0, x)=\phi(x)   \;\;\quad \; & \mbox{on}\;\; \partial \Omega.\end{array}\right.\end{eqnarray}
    ${\mathcal{K}}_m (t,x,y)$ is said to be the heat kernel associated to the perturbed polyharmonic Dirichlet-to-Neumann operator $\Lambda_m$, if for any fixed $y\in \partial \Omega$,
  \begin{eqnarray} \label{55-2} \left\{ \begin{array}{ll} \frac{\partial {\mathcal{K}}_m(t, x,y)}{\partial t} = \Lambda_m {\mathcal{K}}_m(t, x, y),
  \quad  t>0, \,\, x\in \partial \Omega,\\
  {\mathcal{K}}_m (0, x, y)= \delta(x-y)I_m,\end{array}\right. \end{eqnarray}
  where $\Lambda_m$ acts on the $x$ variable and $\delta (x-y)$ is the delta function concentrated at $y$, $I_m$ is the $m\times m$ identity matrix
   (${\mathcal{K}}_m (t,x,y)$ is also called the kernel of the semigroup $e^{t\Lambda_m}$).
    It is well-known that there exists a matrix-valued heat kernel ${\mathcal{K}}_m(t, x, y)$ such that
 \begin{eqnarray}   e^{t\Lambda_m} \phi(x) = \int_{\partial \Omega} {\mathcal{K}}_m (t,x,y) \phi (y) dy, \quad \mbox{for any} \;\; \phi\in (C^\infty(\partial \Omega))^m.\end{eqnarray}
 As be showed in \cite{Gr}, the trace of the associated heat kernel admits an asymptotic expansion
 \begin{eqnarray} \label{1,,1} &&\sum_{k=0}^\infty e^{-t \lambda_k I_m} = Tr\, e^{t\Lambda_m}  \sim \sum_{l=0}^\infty \bigg(\int_{\partial \Omega} a_{m,l}(x) dx\bigg)t^{-n +l}\\
   && \quad \quad \quad \quad \quad \;\quad +\sum_{l=1}^\infty \bigg(\int_{\partial \Omega} b_{m,l} (x)dx\bigg) t^{l} \log t.\nonumber\end{eqnarray}  The matrix-valued coefficients $\int_{\partial \Omega}a_{m,l}(x)dx$ and $\int_{\partial \Omega} b_{m,l}(x)dx$ are called the perturbed polyharmonic Steklov heat invariants, and it follows from (\ref{1,,1})
that they are determined by the perturbed polyharmonic Steklov spectrum.

Heat invariants are metric invariants of the given manifold  $\partial \Omega$ and contain a
lot of information about its geometry and topology (see \cite{PS}, \cite{Liu3}, \cite{Kac}, \cite{Gil}, \cite{EW}, etc.). Higher coefficients are of a considerable interest to physicists since they are connected with many notions of quantum
gravity (see, for example, \cite{Fu}, \cite{ANPS}).
Computation of heat invariants is also a challenging problem in spectral geometry (see \cite{Gr}, \cite{Ta2}, \cite{Gil}, \cite{Liu2}, \cite{PS}, \cite{DH}). This problem is quite similar to the well-known Kac problem for the Laplacian on a domain (The Kac question asks:
is it possible to hear the shape of a domain just by ¡°hearing¡±
all of the eigenvalues of the Dirichlet Laplacian? see \cite{BB}, \cite{vdB},
\cite{R.C}, \cite{Iv}, \cite{Kac}, \cite{LFP}, \cite{AGMT},  \cite{Lo}, \cite{ANPS}, \cite{We1}, \cite{See}, \cite{Sar}, \cite{CH}, \cite{Ch1}, \cite{MS}, \cite{CLN} and the references therein).
 However, the complexity of explicit formulas for $a_{m,l}(x)$ is increasing very rapidly with the growth of $l$.
 For $m=1$ and $q\equiv 0$, Polterovich and Sher \cite{PS}  calculated the first three coefficients $a_{1,l}(x)$ ($0\le l\le 2$), and the author of the present paper \cite{Liu3} obtained the first four coefficients $a_{1,l}(x)$ ($0\le l\le 3$) by a completely different method.

In this paper, we establish a procedure, from which all coefficients $a_{m,l}(x)$ of the asymptotic expansion of the trace of the heat kernel associated to the perturbed polyharmonic Dirichlet-to-Neumann operator can be calculated as $t\to 0^+$. In particular, we also explicitly obtain the first four coefficients $a_{m,l}(x,n)$ ($m\ge 1, \,0\le l\le 3$) of the asymptotic expansion. This provides the information for the area and curvatures of the boundary $\partial \Omega$ by the spectrum of perturbed polyharmonic Steklov problem.

  The main ideas are as follows:
  we first rewrite the perturbed polyharmonic operator as a second order system. By a method of factorization we calculate the full symbol of $\Lambda_m$.
  This method comes from Krupchyk, Lassas and Uhlmann \cite{KLU} (also see, an earlier paper \cite{LU} of Lee and Uhlmann), in which they calculated
  the full symbol in the Euclidean domain. Next, we calculate
   the full symbol of the resolvent operator $(-\Lambda_m-\tau I_m)^{-1}$. However, the computational complexity is amazing high with the growth of $l$.
  In order to overcome this difficult, we decompose the symbol of $(-\Lambda_m-\tau I_m)^{-1}$ into two components: one part depends  on the
  matrix $J_m$ (see Section 3), another part is independent of $J_m$. The first part can be directly calculated, and the second part can be obtained from the author's result in \cite{Liu3}. Finally, by symbol approximation for the formula $e^{t\Lambda_m} =\frac{i}{2\pi} \int_\Gamma (-\Lambda_m-\tau I_m)^{-1} e^{-t\tau} d\tau$, we eventually obtain the following asymptotic expansion:
    \begin{eqnarray*} &&\int_{\partial \Omega} {\mathcal{K}}_m(t,x,x)dS(x) = t^{-n} \int_{\partial \Omega}a_{m,0}(x,n)\, dS(x)+ t^{1-n} \int_{\partial \Omega} a_{m,1}(n,x) \, dS(x) \\
   &&\quad\; + \cdots + t^{M-1-n} \int_{\partial \Omega} a_{m,M-1}(n,x)\, dS(x) +\left\{
           \begin{array}{ll} O(t^{M-n})I_m, \;\; n>M-1\\ O(t\log t)I_m, \;\; n=M-1\end{array}\right. \quad \;\mbox{as}\;\; t\to 0^+,\end{eqnarray*}
            where  $M$ can be taken as $2,3,4,\cdots$ with $n\ge M-1$, and the notation $O(t^j)$ (respectively, $O(t\log t)$) denotes a function which satisfies $|O(t^j)|\le c_0 t^j$ (respectively, $|O(t\log t)|\le c_0t\log t$) for some constant $c_0>0$ and all $t>0$.
           The first coefficient matrix  $a_{m,0}(n,x)=\frac{\Gamma(\frac{n+1}{2})}{\pi^{\frac{n+1}{2}}}I_m$ is independent of $x$. The second coefficient matrix $a_1(n,x)$ depends only on the ``area'' $\mbox{vol}(\partial \Omega)$  and the mean curvature of the boundary $\partial \Omega$. The coefficient matrix $a_2(n,x)$ depends not only on the ``area'' $\mbox{vol}(\partial \Omega)$ and the principal curvatures $\kappa_1,\cdots,\kappa_n$, but also on the scalar curvature ${\tilde{R}}_{\Omega}$ (respectively $R_{\partial \Omega}$) of $\Omega$ (respectively $\partial \Omega$). Finally the fourth coefficient matrix $a_3(n,x)$ depends on $\mbox{vol}(\partial \Omega)$, the principal curvatures, the Ricci tensor ${\tilde{R}}_{jj}$ (respectively $R_{jj}$) and the scalar curvature ${\tilde R}_\Omega$ (respectively $R_{\partial \Omega}$ of $\Omega$ (respectively $\partial \Omega$) as well as the covariant derivative $\sum_{j=1}^n {\tilde R}_{j(n+1)j(n+1),(n+1)}$ of the Ricci curvature with respect to $(\bar \Omega, g)$ (see Theorem 3.1). Generally, $a_{m,M-1} (n,x)$ can all be represented in terms of the geometical quantities (see, (ii) of Remark 3.2).
             This asymptotic expansion shows that one can hear the ``area'' of $\partial \Omega$ and all $\int_{\partial \Omega} a_{m,M-1}(n,x)dS(x)$ ($M=1,2,3,\cdots$)  by ``hearing'' all of the Steklov eigenvalues.
            Since $\int_{\partial \Omega} a_{m,M-1}(n,x) dS(x)$ ($M\ge 1$) are all spectral invariants, it immediately  follows that two domains with different spectral invariants can never have the same perturbed polyharmonic Steklov spectrum.

\vskip 1.49 true cm

\section{Full symbol of the perturbed polyharmonic Dirichlet-to-Neumann operator }

\vskip 0.45 true cm

We  rewrite the equation (\ref{1.1}) as a second order system. Introducing
$$ u_1= u, \; u_2=
(-\Delta_g)u,\, \cdots,\, u_m=(-\Delta_g)^{m-1} u,$$
we get
\begin{eqnarray} \label {2-1}  \big(-\Delta_g \otimes I_m + J_m\big) U=0\quad \; \mbox{in}\;\; \Omega,\end{eqnarray}
where $U=(u_1, u_2, \cdots, u_m)^t$ and
 \begin{gather*} \label{3--1}
   J_m=\begin{pmatrix} 0 & -1 & 0 &\cdots   &0\\
     0 &  0 & -1 &\cdots  &  0\\
      \vdots & \vdots & \vdots & \ddots & \vdots  \\
     0 & 0 & 0& \cdots &  -1\\
      q(x) & 0 &0&  \cdots   &   0  \end{pmatrix}. \end{gather*}
 The graph of the above system (\ref{2-1}) is  defined as
$$\{ (U\big|_{\partial \Omega}, \partial_\nu U\big|_{\partial \Omega})\big|\, U\in \prod_{j=1}^{m-1} H^{m-2j-1/2} (\partial \Omega), \;\; U \;\;\mbox{solves}\;\; (\ref{2-1})\}.$$
Obviously,  the graph of the
system (\ref{2-1}) and the graph of the equation (\ref{1.3*}) coincide.

We will compute full symbol of the perturbed polyharmonic Dirichlet-to-Neumann operator $\Lambda_m$ by following the method of \cite{KLU}.
When considering the system (\ref{2-1}) near the boundary, we  make use of the boundary  normal coordinates.
For each $z\in \partial \Omega$, let $x_{n+1}$ denote the parameter along the unit-speed geodesic starting at $z$
  with initial direction given by the inward boundary normal to $\partial \Omega$. If $x=(x_1, \cdots, x_n)$ are any local coordinates for
 $\partial \Omega$ near $z\in \partial \Omega$, we can extend them smoothly to functions on a neighborhood
of $z$ in $\Omega$ by letting them be constant along each normal geodesic.  It follows easily that $(x, x_{n+1})$
form coordinates for $\bar \Omega$ in some neighborhood of $z$, which we call the boundary
normal coordinates determined by $x$. In these coordinates $x_{n+1}>0$ in
 $\Omega$, and $\partial \Omega$ is locally characterized by $x_{n+1}=0$ (see, \cite{LU} or \cite{Ta2}).
    Then the metric tensor on $\bar \Omega$ has
the form (see, p.$\,$532 of \cite{Ta2})
\begin{gather} \label{a-1} \big(g_{jk} (x,x_{n+1}) \big)_{(n+1)\times (n+1)} =\begin{pmatrix} ( h_{jk} (x,x_{n+1}))_{n\times n}& 0\\
      0& 1 \end{pmatrix}. \end{gather}

      In what follows, we will use the notations $\partial_{x_j}=\partial/\partial x_j, \, D_{x_j}= -i \partial_{x_j}$, and $D =(D_{x_1}$, $\cdots$, $D_{x_n}$, $D_{x_{n+1}})$. The Laplace-Beltrami operator $\Delta_g$ on $\Omega$ is given in local coordinates by
\begin{eqnarray*} \label{a-2} \Delta_g u&=&\sum_{j,k=1}^{n+1} \frac{1}{\sqrt{|g|}}\, \frac{\partial}{\partial x_j}\left( \sqrt{|g|}\, g^{jk} \frac{\partial u}{\partial x_k}\right) \\
&=& \frac{\partial^2 u} {\partial x_{n+1}^2} +\frac{1}{2} \frac{\partial \log |h|}{\partial x_{n+1} } \, \frac{\partial u}{\partial x_{n+1}}\\
&&  + \sum_{j,k=1}^n \left( h^{jk}\frac{\partial^2 u} {\partial x_j\partial x_j} +\frac{1}{2} h^{jk} \frac{\partial \log|h|}{\partial x_j}
\, \frac{\partial u}{\partial x_k}  +\frac{\partial h^{jk}}{\partial x_j} \, \frac{\partial u}{\partial x_k}\right). \end{eqnarray*}
So one can write (see, p.$\,$1101 of \cite{LU})
  $$ -\Delta_g = D_{x_{n+1}}^2 +i E(x,x_{n+1}) D_{x_{n+1}} + Q(x, x_{n+1},D_{x}),$$
where \begin{eqnarray*}   E(x,x_{n+1}) &=&  -\frac{1}{2} \sum_{j,k=1}^n h^{jk}(x,x_{n+1}) \, \frac{\partial  h_{jk}(x,x_{n+1})}{\partial x_{n+1}},\\
Q(x, x_{n+1}, D_{x}) &=& \sum_{j,k=1}^n h^{jk} (x,x_{n+1}) D_{x_j} D_{x_k} \\
  && -i \sum_{j,k=1}^n  \left( \frac{1}{2} h^{jk} (x,x_{n+1}) \, \frac{\partial \log|h(x,x_{n+1})|}{\partial x_j} +\frac{\partial h^{jk} (x,x_{n+1})}{\partial x_j} \right) D_{x_k}.\end{eqnarray*}
 The principal part
$Q_2$ of $Q$ is given by
$$ Q_2(x,x_{n+1}, D_{x¡¯}) =\sum_{j,k=1}^{n} h^{jk} (x,x_{n+1}) D_{x_j} D_{x_k}.$$
Here $(h^{jk})$  is the inverse of the matrix $(h_{jk})$.
In the boundary normal coordinates, the operator in the system (\ref{2-1}) has the form
\begin{eqnarray}\label{2-4} && P(x,x_{n+1},D) := \big(D_{x_{n+1}}^2 + i E(x,x_{n+1}) D_{x_{n+1}} +Q(x,x_{n+1}, D_{x})\big)\otimes I_m   +J_m. \end{eqnarray}

  The following result is due to Krupchyk, Lassas and Uhlmann (see, Proposion 4.2 of \cite{KLU}).

\vskip 0.2 true cm

\noindent  {\bf Lemma 2.1.} \    {\it There is a matrix-valued pseudodifferential operator $B(x, x_{n+1}, D_{x})$
of order one in $x$ depending smoothly on $x_{n+1}$ such that
\begin{eqnarray}\label{2-5}  && P(x,x_{n+1}, D)=\big(D_{x_{n+1}} \otimes I_m +iE(x,x_{n+1})\otimes I_m    -iB(x,x_{n+1}, D_{x})\big)\\
&&\quad \;\quad \quad \qquad \qquad \, \times  \big(D_{x_{n+1}}\otimes I_m +iB(x, x_{n+1}, D_{x})\big),\nonumber\end{eqnarray}
modulo a smoothing operator. Here $B(x,x_{n+1}, D_{x})$ is unique modulo a smoothing term,
if we require that its principal symbol is given by
$-\sqrt{Q_2(x,x_{n+1},\xi)}\; I_m$}.

\vskip 0.2 true cm

By combining (\ref{2-4}) and (\ref{2-5}), we have
\begin{eqnarray}\label{2-6} && \quad B^2 (x,x_{n+1}, D_{x}) +i
[D_{x_{n+1}} \otimes I_m, B(x,x_{n+1}, D_{x})] - E(x,x_{n+1}) B(x, x_{n+1},D_{x})\\
 && \quad \quad \quad \qquad   = Q(x,x_{n+1}, D_{x})\otimes I_m  +J_m,\nonumber \end{eqnarray}
modulo a smoothing operator. Let us write the full symbol of $B(x, x_{n+1},D_{x})$ as
$$B(x, x_{n+1}, \xi) \sim \sum_{j=-\infty}^1 {\tilde{r}}_j(x,x_{n+1}, \xi),$$
with ${\tilde{r}}_j$ taking values in $m\times m$ matrices with entries homogeneous of degree $j$ in
$\xi=(\xi_1,\cdots,\xi_{n})\in {\Bbb R}^{n}$. Thus, (\ref{2-6}) implies (see \cite{KLU}) that
\begin{eqnarray} \label{2-7}  &&  \sum_{l=-\infty}^2 \big( \sum_{\underset {|\alpha|\ge 0, \, j,k\le 1}  {j+k-|\alpha|=l}}
 \frac{1}{\alpha!} \partial^\alpha_{\xi} {\tilde{r}}_j  D^\alpha_{x} {\tilde{r}}_k \big) -E(x,x_{n+1}) \sum_{j=-\infty}^1 {\tilde{r}}_j(x,x_{n+1}, \xi) \\ && \qquad  \quad+\sum_{j=-\infty}^1 \partial_{x_{n+1}} {\tilde{r}}_j (x, x_{n+1},\xi)    = Q_2(x,x_{n+1}, \xi)I_m + Q_1 (x,x_{n+1},  \xi)I_m  +J_m.\nonumber \end{eqnarray}
Here $Q_1(x, x_{n+1},\xi)=Q(x,x_{n+1}, \xi)-Q_2(x,x_{n+1},\xi)$
 is homogeneous of degree one in $\xi$. Equating
the terms homogeneous of degree two in (\ref{2-7}), we get
$$ {\tilde{r}}_1^2(x,x_{n+1}, \xi) = Q_2(x,x_{n+1}, \xi)I_m ,$$
so we may choose ${\tilde{r}}_1(x, x_{n+1}, \xi)$ such that (see \cite{KLU})
\begin{eqnarray}\label{2-8-} {\tilde{r}}_1(x,x_{n+1}, \xi) =-\sqrt{Q_2(x, x_{n+1}, \xi)}\;I_m.\end{eqnarray}
Equating the terms homogeneous of degree one in (\ref{2-7}), we have an equation
\begin{eqnarray} \label{2-8}&&
2{\tilde{r}}_1 {\tilde{r}}_0+ \sum_{|\alpha|=1} \partial^\alpha_{\xi} {\tilde{r}}_1 \cdot D^\alpha_{x} {\tilde{r}}_1 -E{\tilde{r}}_1 +\partial_{x_{n+1}} {\tilde{r}}_1 = Q_1 (x,x_{n+1}, \xi)I_m.\nonumber\end{eqnarray}
By considering the terms homogeneous of degree zero in (\ref{2-7}), we get
\begin{eqnarray} \label{2-9}  && {\tilde{r}}_0^2 +2 {\tilde{r}}_{-1} {\tilde{r}}_1 +\sum_{|\alpha|=1} \big(\partial^\alpha_{\xi} {\tilde{r}}_1 \cdot D_{x}^\alpha {\tilde{r}}_0 + \partial^\alpha_{\xi} {\tilde{r}}_0\cdot D^\alpha_{x} {\tilde{r}}_1\big) +\sum_{|\alpha|=2} \frac{1}{\alpha!} \partial^\alpha_\xi {\tilde{r}}_1 \cdot D_x^\alpha {\tilde{r}}_1 \\ && \quad \quad \quad
 -E {\tilde{r}}_0 +\partial_{x_{n+1}} {\tilde{r}}_0 =J_m. \nonumber\end{eqnarray}
Also, from the  terms homogeneous of degree $-1$ in (\ref{2-7}), it follows that
\begin{eqnarray}\label{2v10}  && \quad  2({\tilde{r}}_{-2} {\tilde{r}}_{1} +{\tilde{r}}_{-1} {\tilde{r}}_0 ) +\sum_{|\alpha|=1} \big(\partial^\alpha_{\xi} {\tilde{r}}_0  \cdot D^\alpha_{x}{\tilde{r}}_0
 +\partial_\xi^\alpha {\tilde{r}}_1 \cdot D_x^\alpha {\tilde{r}}_{-1} +\partial_\xi^\alpha {\tilde{r}}_{-1} \cdot D_x^\alpha {\tilde{r}}_1\big)\\
  && \quad \quad \quad \quad +
    \sum_{|\alpha|=2} \frac{1}{\alpha!}\big(\partial^\alpha_{\xi} {\tilde{r}}_0 \cdot  D^\alpha_{x}{\tilde{r}}_1 + \partial_\xi^\alpha {\tilde{r}}_1 \cdot D_x^\alpha {\tilde{r}}_0\big)
    +\sum_{|\alpha|=3} \frac{1}{\alpha!}\, \partial^\alpha_\xi {\tilde{r}}_1 \cdot D_x^\alpha {\tilde{r}}_1\nonumber\\ &&
  \quad \quad \quad \quad  -E {\tilde{r}}_{-1} +\partial_{x_{n+1}} {\tilde{r}}_{-1}=0.\nonumber \end{eqnarray}
The terms ${\tilde{r}}_j, \, j\le -3$, are chosen in a similar fashion,
by equating terms of degree of homogeneity $j+1$ in (\ref{2-7}) (cf. \cite{KLU}).

Clearly, the perturbed polyharmonic Dirichlet-to-Neumann operator $\Lambda_m$ is just the restriction of $B(x,x_{n+1},D_x)$ to $x_{n+1}=0$. Write the symbol of the operator $-\Lambda_m$ as $r_1+r_0+r_{-1}+\cdots $.
    From (\ref{2-8-})---(\ref{2v10}), we obtain
      \begin{eqnarray*} && r_1= \sqrt{Q_2(x,0,\xi)}\; I_m,\\
         && r_0=-\frac{1}{2} r_1^{-1}\bigg[ Q_1(x,0,\xi)I_m -\sum_{|\alpha|=1} \partial_{\xi}^\alpha r_1 \cdot D^\alpha_{x} r_1 +E r_1 -\partial_{x_{n+1}} r_1 \bigg],\end{eqnarray*}
 \begin{eqnarray*}  &&r_{-1}=\frac{1}{2} r_1^{-1} \bigg[- J_m +r_0^2 +\sum_{|\alpha|=1} \big(\partial_{\xi}^\alpha r_1 \cdot  D^\alpha_{x} r_0
 + \partial^\alpha_{\xi} r_0 \cdot D^\alpha_{x} r_1\big)\quad \quad  \\
 && \quad \quad \quad +\sum_{|\alpha|=2} \frac{1}{\alpha!} \partial_\xi^\alpha r_1 \cdot D_x^\alpha r_1  -E r_0 +\partial_{x_{n+1}} r_0 \bigg],\nonumber \end{eqnarray*}
    \begin{eqnarray*} &&  r_{-2}=\frac{1}{2}r_1^{-1} \bigg[2r_{-1}r_0 +\sum_{|\alpha|=1}\big( \partial_{\xi}^\alpha r_0 \cdot D_x^\alpha r_0 +\partial_\xi^\alpha r_1 \cdot D_x^\alpha r_{-1}  + \partial_\xi^\alpha r_{-1} \cdot D_x^\alpha r_1\big) \\
    && \; \quad \quad \; +
 \sum_{|\alpha|=2} \frac{1}{\alpha!} \big(\partial^\alpha_{\xi} r_0 \cdot  D^\alpha_{x}r_1 +
 \partial^\alpha_{\xi} r_1 \cdot  D^\alpha_{x}r_0\big) \\
  && \; \quad \quad\; +\sum_{|\alpha|=3} \frac{1}{\alpha!}\, \partial^\alpha_\xi r_1\cdot D_x^\alpha r_1  -E r_{-1} +\partial_{x_{n+1}} r_{-1}\bigg].\nonumber \end{eqnarray*}

\vskip 0.25 true cm

For the Laplacian $\Delta_g$ on a  compact  Riemannian manifold without boundary, the first method for derivation of heat kernel asymptotics is due to Seeley (\cite{See}). This method was developed later by Gilkey (see, Theorem 1.3 in \cite{Gil}) who
presented a way to get recursive formulas for the heat invariants.
We now apply the method of Seeley \cite{See} (cf. the paper of Gilkey and Grubb \cite{GiGr}) to calculate the asymptotic expansion of the trace of the heat kernel for the perturbed polyharmonic Dirichlet-to-Neumann operator $\Lambda_m$.
Let $\epsilon >0$ be given. The spectrum of $-\Lambda_m$ lies in a cone of slope $\epsilon$ about the positive real axis. Let
 $\Gamma$ be a path about the cone with slope $2\epsilon$ outside some compact set. For $\tau$
on $\Gamma$, the operator $(-\Lambda_m-\tau I_m)^{-1}$ is a uniformly bounded compact operator from
 $(L^2(\partial \Omega)^m\to (L^2(\partial \Omega))^m$. The integral
  \begin{eqnarray} \label{3>1} \frac{i}{2\pi} \int_{\Gamma} e^{-t\tau I_m} (-\Lambda_m - \tau I_m)^{-1} d\tau\end{eqnarray}
converges absolutely for $t>0$ and defines the operator $e^{t\Lambda_m}$.
We construct a matrix-valued pseudodifferential operator to approximate the resolvant $(-\Lambda_m -\tau I_m)^{-1}$ as follows: let  $s(\tau)$ be a matrix-valued pseudodifferential operator of order $-1$ with parameter $\tau$ for which
\begin{eqnarray*} (-\Lambda_m-\tau I_m) s(\tau) = s(\tau)(-\Lambda_m-\tau I_m) =I_m.\end{eqnarray*}
(Actually, we require that $(-\Lambda_m-\tau I_m) s(\tau)-I_m$ and $s(\tau)(-\Lambda_m-\tau I_m) -I_m$ are both pseudodifferential  operator of order $-\infty$.)  It follows from \cite{GiGr} (see also, \cite{Gil}, \cite{Gr} or \cite{See}) that such a pseudodifferential operator must have symbol $s_{-1} (x, \xi,\tau) +s_{-2} (x, \xi, \tau) +s_{-3}  (x, \xi, \tau) +\cdots$
given by \begin{eqnarray*}
 & s_{-1} (x, \xi, \tau) = (r_1 -\tau I_m)^{-1},
\\ & s_{-1-l} (x, \xi, \tau) =- (r_1- \tau I_m)^{-1} \left(\sum_{\underset {|\alpha|= l+j+k\ge 0}{-l \le k\le 1,\; -l\le j \le -1}}
 \frac{1}{\alpha!} \, \partial^\alpha_\xi r_k \cdot D_x^\alpha s_j\right).\end{eqnarray*}
 For the sake of convenience, we write out the expressions for $s_{-2}, s_{-3}$ and $s_{-4}$:
 \begin{eqnarray*}&& s_{-2} = - (r_1 -\tau I_m)^{-1}r_0 s_{-1}- (r_1-\tau I_m)^{-1} \sum_{|\alpha|=1} \partial^\alpha_\xi r_1 \cdot D_x^\alpha s_{-1}, \\
    && s_{-3} =-(r_1-\tau I_m)^{-1} \bigg[ r_0 s_{-2} +r_{-1} s_{-1} +\sum_{|\alpha|=1} \big(\partial_\xi^\alpha r_1 \cdot D_x^\alpha s_{-2} + \partial_\xi^\alpha r_0\cdot D_x^\alpha s_{-1}\big)  \\
    && \quad \quad \quad  +\sum_{|\alpha|=2} \frac{1}{\alpha!} \partial^\alpha_\xi r_1  \cdot D_x^\alpha s_{-1}\bigg],
     \\  && s_{-4} = -(r_1-\tau I_m)^{-1}  \bigg[ r_0 s_{-3}  +r_{-1} s_{-2} +r_{-2} s_{-1} \\
     && \quad \quad \quad +\sum_{|\alpha|=1} \big(\partial_\xi^\alpha r_1  \cdot D_x^\alpha s_{-3} +\partial_\xi^\alpha r_0  \cdot D_x^\alpha s_{-2} + \partial_\xi^\alpha r_{-1} \cdot D_x^\alpha s_{-1}\big) \\
      && \quad \quad \quad + \sum_{|\alpha|= 2} \frac{1}{\alpha!} \big(\partial_\xi^\alpha r_1  \cdot D_x^\alpha s_{-2} + \partial^\alpha_\xi r_0  \cdot D_x^\alpha s_{-1} \big) +\sum_{|\alpha|=3}  \frac{1}{\alpha!} \partial_\xi^\alpha r_1 \cdot D_x^\alpha s_{-1}\bigg].\end{eqnarray*}
Let  $T^*_x(\partial \Omega)$ be the cotangent space at $x$. It follows from (\ref{3>1}) that for $m\ge 1,\; 0\le l\le n$,
 \begin{eqnarray} \label{3./1} && a_{m,l}(x) =\frac{i}{(2\pi)^{n+1}}  \int_{T^*_x (\partial \Omega)} \int_\Gamma e^{-\tau I_m} s_{-1-l} (x, \xi,\tau) d\tau \, d\xi\\
&& \quad \quad \;\; =\frac{1}{(2\pi)^{n}}  \int_{{\Bbb R}^n} \bigg(\frac{i}{2\pi} \int_\Gamma e^{-\tau } s_{-1-l} (x, \xi,\tau) d\tau \bigg) d\xi .\nonumber\end{eqnarray}

\vskip 1.49 true cm

\section{Calculation of heat invariants}

\vskip 0.45 true cm

 In order to state our main result regarding the relationship between the spectrum of the perturbed polyharmonic Dircichlet-to-Neumann  operator and various geometric quantities on $\partial \Omega$, we introduce the following notations:
 we denote by $\,{\tilde R}_{jkjk}(x)$ (respectively $R_{jkjk}(x)$), $\,{\tilde{R}}_{jj}(x)=\sum_{k=1}^{n+1} {\tilde{R}}_{jkjk}(x)$ (respectively
 ${R}_{jj}(x)=\sum_{k=1}^{n} R_{jkjk}(x)$), $\, {\tilde {R}}_{\Omega}$ (respectively $R_{\partial \Omega} (x)$)
 the curvature tensor, the Ricci curvature tensor, the scalar curvature with respect to $\Omega$ (respectively $\partial \Omega$) at $x\in \partial \Omega$. Denote by  $\,\sum_{j=1}^n {\tilde R}_{j(n+1)j(n+1), (n+1)} (x)$ the covariant derivative of the curvature tensor with respect to $(\bar \Omega, g)$ in the inward normal direction $\nu$. In addition, we denote by  $\kappa_1 (x), \cdots, \kappa_n(x)$  the principal curvatures of $\partial \Omega$ at $x\in \partial \Omega$,
and $dS(x) =\sqrt{|h(x)|}dx$ the area element of $\partial \Omega$. Recall also that \begin{gather*}
J_m(x)=\begin{pmatrix} 0 & -I_{m-1}\\
      q(x)& 0\end{pmatrix}. \end{gather*}

\vskip 0.29 true cm

\noindent{\bf Theorem 3.1.} \ {\it Suppose that $(\mathcal{M},g)$ is an $(n+1)$-dimensional, smooth Riemannian
manifold, and assume that $\Omega\subset \mathcal{M}$ is a bounded domain
with smooth boundary $\partial \Omega$.  Let ${\mathcal {K}}_m(t, x, y)$ be the heat kernel associated to the perturbed polyharmonic Dirichlet-to-Neumann operator
$\Lambda_m$ ($m\ge 1$).

   (a)  \  If $n\ge 1$, then
   \begin{eqnarray} \label{6.0.1}    && \int_{\partial \Omega} {\mathcal{K}}_m (t, x,x)dx=t^{-n}
\int_{\partial \Omega} a_{m,0}(n,x) \,dS(x) + t^{1-n}  \int_{\partial \Omega} a_{m,1}(n,x)\, dS(x) \\  &&\quad \quad \quad \quad + \left\{\begin{array}{ll} O(t^{2-n}) \quad \;\, \mbox{when}\;\; n>1,\\ O(t\log t) \quad \, \mbox{when} \;\; n=1,\end{array}\right. \quad \;\;\mbox{as}\;\; t\to 0^+;\nonumber\end{eqnarray}

   (b) \  If $n\ge 2$, then \begin{eqnarray} \label{6.0.1-2}    && \int_{\partial \Omega} {\mathcal{K}}_m (t, x,x)dx=t^{-n}
\int_{\partial \Omega} a_{m,0}(n,x) \,dS(x) + t^{1-n}  \int_{\partial \Omega} a_{m,1}(n,x)\, dS(x)\\
  && + t^{2-n} \int_{\partial \Omega} a_{m,2}(n,x) \,dS(x) + \left\{\begin{array}{ll} O(t^{3-n}) \quad \, \mbox{when}\;\; n>2,\\ O(t\log t) \; \;\;\mbox{when} \;\; n=2,\end{array}\right.    \quad \;\;\mbox{as}\;\; t\to 0^+;\nonumber\end{eqnarray}

    (c) \  If $n\ge 3$, then
\begin{eqnarray} \label{6.0.1-3}    && \int_{\partial \Omega} {\mathcal{K}}_m (t, x,x)dx=t^{-n}
\int_{\partial \Omega} a_{m,0}(n,x) \,dS(x) + t^{1-n}  \int_{\partial \Omega} a_{m,1}(n,x)\, dS(x) \\
  && \qquad \qquad  + t^{2-n} \int_{\partial \Omega} a_{m,2}(n,x) \,dS(x)  + t^{3-n} \int_{\partial \Omega} a_{m,3}(n,x) \,dS(x) \nonumber\\&&\qquad \qquad +
  \left\{\begin{array}{ll} O(t^{4-n}) \quad\, \, \mbox{when}\;\; n>3,\\ O(t\log t) \quad \mbox{when} \;\; n=3,\end{array}\right.
   \quad \;\;\mbox{as}\;\; t\to 0^+.\nonumber\end{eqnarray} Here
   \begin{eqnarray} \label{06-060-1} a_{m,0}(n,x)=\frac{\Gamma(\frac{n+1}{2})}{\pi^{\frac{n+1}{2}}}I_m,\end{eqnarray}
  \begin{eqnarray}\label{06-060}\quad a_{m,1}(n,x)= \left(\frac{1}{2\pi}\right)^{n}   \frac{(n-1)\Gamma(n)\,\mbox{vol}({\Bbb S}^{n-1})}{2n} \bigg(\sum_{j=1}^n \kappa_j(x)\bigg)I_m,
     \end{eqnarray}
    \begin{eqnarray}  \label{66.55}&&    a_{m,2}(n,x) = \frac{\Gamma(n-1) \cdot vol({\Bbb S}^{n-1})}{2(2\pi)^n}\, J_m \\
   &&\quad \qquad \qquad \;  + \frac{\Gamma(n-1) \mbox{vol}({\Bbb S}^{n-1}) }{8(2\pi)^n} \left[ \frac{3-n}{3n} R_{\partial \Omega} + \frac{n-1}{n} {\tilde R}_{\Omega}\right. \nonumber\\ && \left.\qquad \qquad \quad \;  +\frac{n^3 -n^2 -4n +6}{n(n+2)} \big(\sum_{j=1}^n \kappa_j (x)\big)^2 + \frac{n^2 -n-2}{n(n+2)} \sum_{j=1}^n \kappa_j^2(x) \right]I_m
      \nonumber \end{eqnarray} and
\begin{eqnarray} \label{6.6-1.} && a_{m,3}(n,x)= \left(\frac{1}{2\pi}\right)^n \bigg(\frac{n-1}{4n} \Gamma(n-1)   -\frac{1}{2} \Gamma(n-2) \bigg) vol({\Bbb S}^{n-1})  \bigg(\sum_{j=1}^n  \kappa_j\bigg)J_m. \\
&&  \quad \, + \left(\frac{1}{2\pi}\right)^n \frac{\Gamma(n-2)}{4} vol({\Bbb S}^{n-1}) \big(\partial_{x_{n+1}} J_m\big)\nonumber\\
&&  \quad \, +  \bigg(\frac{1}{2\pi}\bigg)^{n}\frac{\Gamma(n-2)\,\mbox{vol} (S^{n-1}) }{8n}  \bigg[
   \frac{n^3 -2n^2 -7n+7}{2(n+2)}  {\tilde R}_\Omega(x)\, \big(\sum_{j=1}^n \kappa_j(x)\big) \nonumber \\ && \quad \, + \frac{-3n^4 -4 n^3 +59n^2 +75 n -180}{ 6(n+2)(n+4)} \big(\sum_{j=1}^n \kappa_j (x) \big)R_{\partial \Omega} \nonumber\\ && \quad \,+
   \frac{n^5 -20 n^3 +2n^2 +61 n -74}{6(n+2)(n+4)} \big(\sum_{j=1}^n \kappa_j(x)\big)^3 \nonumber\\ && \quad \,
    + \frac{n^4 +8n^3 +15n^2 +3n -32}{2(n+2) (n+4)} \big(\sum_{j=1}^n \kappa_j(x)\big)\big(\sum_{j=1}^n \kappa_j^2 (x)\big)\nonumber\\ && \quad \,
    +\frac{-6n^3-34n^2 +40}{3(n+2)(n+4)} \sum_{j=1}^n \kappa_j^3(x) + \frac{4n^2-6}{n+2} \sum_{j=1}^n\kappa_j(x) {\tilde R}_{jj}(x)\nonumber
    \\ && \quad \,-\frac{12n^3 +50 n^2 -6n -104}{3(n+2)(n+4)} \sum_{j=1}^n \kappa_j(x) R_{jj}(x) + (n-1) \sum_{j=1}^n {\tilde{R}}_{j(n+1)j(n+1),(n+1)}(x)\nonumber\\ &&\quad \,-\frac{n-2}{2 } {\tilde{R}}_{\Omega} (x)  +\frac{n-2}{2} R_{\partial \Omega} -\frac{n-2}{2} \sum_{j=1}^n \kappa_j^2(x)\bigg]I_m. \nonumber \end{eqnarray}}

\noindent  {\it Proof.} \ Since  $Z=\int_{\partial \Omega} {\mathcal{K}}_m(t, x, x)dS(x)$ converges, $e^{t \Lambda_m} : \phi \to \int_{\partial \Omega}
 {\mathcal{K}}_m\phi$ is a compact mapping of the (real) Hilbert space $H=(L^2 (\partial \Omega, dS(x)))^m$.
This implies \begin{eqnarray} \label {413} {\mathcal{K}}_m(t,x,y)=\sum_{k=1}^{\infty} e^{-t\lambda_k I_m} (\gamma u_k(x))(\gamma u_k(y))  \end{eqnarray}
with uniform convergence on compact figures of
$(0, \infty)\times \partial \Omega\times \partial \Omega$,
and the spur $Z$ is easily evaluated as (see, for example,   Chapter 4 of \cite{Gr})
\begin{eqnarray}  Z= \int_{\partial \Omega} \sum_{k=1}^\infty e^{-t\lambda_k I_m} (\gamma u_k(x))^2dS(x)
= \sum_{k=1}^\infty e^{-t\lambda_k I_m}.\end{eqnarray}
Thus, our task is  to estimate the pole ${\mathcal{K}}_m (t,x,x)$  for $t\downarrow 0$,  up to needed remainder term.

Let us split $r_j$ ($j\le -1$) into two components $r'_{j}$ and $r''_{j}$, where $r'_{j}$ depends on $J_m$, and $r''_{j}$ is independent of $J_m$  (here $r_j$  are given in Section 2).  We further choose our boundary normal coordinates frame such that the second fundamental form is  diagonal. Thus, for any fixed point $p\in \partial \Omega$ whose coordinate is $(x_0,0)$, by taking $x_0$-centered normal coordinates of $\partial \Omega$ we have  \begin{eqnarray} \label{-4.1} \frac{\partial r_1}{\partial x_l} (x_0,0)=\frac{1}{2} r_1^{-1} \sum_{j,k=1}^n \frac{\partial h^{jk}}{\partial x_l} (x_0,0)=0.\end{eqnarray}
 It then follows that  $\left(-\frac{1}{2}\,\frac{\partial h_{jk}}{\partial x_{n+1}}(x_0,0)\right)_{n\times n}$  is the matrix of Weingarten's map under a basis of $T_0(\partial \Omega)$ at $x_0$, and its eigenvalues are just the principal curvatures $\kappa_1,\cdots,\kappa_n$ of $\partial \Omega$  at $x_0$, i.e.,
 \begin{eqnarray}\label{-4.2} \left( \frac{\partial h_{jk}}{\partial x_{n+1}} (x_0,0)\right)_{n\times n} = \mbox{diag}(-2\kappa_1,\cdots, -2\kappa_n).\end{eqnarray}
  Noting  also that \begin{eqnarray} \label{-4.3} h^{jk} (x_0,0) =\delta^{jk} \quad \; \mbox{and}\;\; \frac{\partial h^{jk}}{\partial x_{n+1}} (x_0, 0) =-\frac{\partial h_{jk}}{\partial x_{n+1}}(x_0,0)=2\kappa_j\delta_{jk}.\end{eqnarray}
       We  find that
  \begin{eqnarray*} && r_1(x_0, \xi)=\sqrt{Q_2(x,0, \xi)}\; I_m =\sqrt{\sum_{j,k=1}^n h^{jk}(x_0) \xi_j\xi_k} \; I_m=r''_1(x_0, \xi),\\
  && r_0(x_0, \xi)=\bigg[ \frac{-i} {8 Q_2^{3/2} }  \sum_{l=1}^n \bigg( \sum_{j,k=1}^n h^{jk} \big(\delta^{lj} \xi_k +\delta^{lk} \xi_j\big) \bigg)\big(\sum_{j,k=1}^n \frac{\partial h^{jk}}{\partial x_l}
 \xi_j\xi_k \big) +\frac{1}{4Q_2} \sum_{j,k=1}^n \frac{\partial h^{jk}}{\partial x_{n+1}} \xi_j\xi_k\\
  && \quad \quad \quad  \quad \quad  + \frac{i}{2\sqrt{Q_2}} \sum_{j,k=1}^{n}\bigg(\frac{1}{2} h^{jk}  \frac{\partial \log |h|}{\partial x_j} +\frac{\partial h^{jk} }{\partial x_j} \bigg) \xi_k +\frac{1}{4} \sum_{j,k=1}^n h^{jk}
  \frac{\partial h_{jk}}{\partial x_{n+1}}\bigg]I_m=r''_0(x_0,\xi),\\
  && r_{-1}(x_0, \xi) =-
  \frac{1}{2} r_1^{-1} J_m + \frac{1}{2} r_1^{-1} \left( r_0^2 +\sum_{|\alpha|=1} \big(\partial^\alpha_\xi r_1 \cdot D^\alpha_x r_0 + \partial_\xi^\alpha r_0\cdot
  D_x^\alpha r_1\big)\right. \\
   && \left.\quad \quad \quad\;\; \qquad +\sum_{|\alpha|=2} \frac{1}{\alpha!} \partial^\alpha_\xi r_1 \cdot D_x^\alpha r_1 -E r_0 +\partial_{x_{n+1}} r_0\right)\\
   && \quad \quad \quad \quad \;\; = -\frac{1}{2} r_1^{-1} J_m +r''_{-1}(x_0, \xi),\end{eqnarray*}
   \begin{eqnarray*}  &&  r_{-2} (x_0, \xi)=  \frac{1}{2} r_1^{-1} \bigg[ 2\big(- \frac{1}{2} r_1^{-1} J_m +r''_{-1}\big) r_0 +\sum_{|\alpha|=1} \bigg( \partial_\xi^\alpha r_0\cdot D_x^\alpha r_0 + \partial_\xi^\alpha r_1 \cdot D_x^\alpha\big(-\frac{1}{2} r_1^{-1} J_m+r''_{-1} \big)\\
   && \qquad \qquad\quad  +\partial_\xi^\alpha \big(-\frac{1}{2} r_1^{-1} J_m +r''_{-1}\big)\cdot D_x^\alpha r_1 \bigg) +\sum_{|\alpha|=2} \frac{1}{\alpha!} \big( \partial_\xi^\alpha r_0 \cdot D_x^\alpha r_1 +\partial_\xi^\alpha r_1 \cdot D_x^\alpha r_0\big)\\
  &&  \qquad \qquad \quad    +\sum_{|\alpha|=3} \frac{1}{\alpha!}  \partial_\xi^\alpha  r_1\cdot D_x^\alpha r_1
   -E \big( -\frac{1}{2} r_1^{-1} J_m +r''_{-1}\big) +\partial_{x_{n+1}} \big(-\frac{1}{2} r_1^{-1} J_m +r''_{-1}\big)\bigg]\end{eqnarray*}
  \begin{eqnarray*}  &&  =\frac{1}{2} r_1^{-1} \bigg[ -  r_1^{-1} J_m r_0 -\frac{1}{2}\sum_{|\alpha|=1}\big( \partial_\xi^\alpha r_1\cdot  r_1^{-1} D_x^\alpha  J_m \big) \qquad \qquad \qquad \qquad \qquad\; \quad \\ && \quad +  \frac{E}{2}r_1^{-1} J_m -\frac{1}{2}\partial_{x_{n+1}} \big( r_1^{-1} J_m\big)\bigg] + r''_{-2}(x_0, \xi).\end{eqnarray*}
  Similarly, we can split $s_j$ $(j\le -1$) into two parts, the first part $s'_j$ depends on $J_m$ and the second part $s''_j$  is independent of $J_m$, i.e.,
  \begin{eqnarray} &&\label{-4.11} s_{-1} =  (r_1-\tau I_m)^{-1} =s''_{-1},\\
   && s_{-2}=  -(r_1-\tau I_m)^{-1} r_0 s_{-1} -(r_1 -\tau I_m)^{-1} \sum_{|\alpha|=1} \partial^\alpha_\xi r_1 \cdot D_x^\alpha s_{-1} =s''_{-2},\nonumber \end{eqnarray}
  \begin{eqnarray*} && s_{-3} = -(r_1-\tau I_m)^{-1} \bigg[ r_0 s_{-2} +\big(-\frac{1}{2} r_1^{-1} J_m +r''_{-1}\big)s_{-1}  +
  \sum_{|\alpha|=1} \big(\partial_\xi^\alpha r_1 \cdot D_x^\alpha s_{-2}\\
   && \quad \quad \quad + \partial_\xi^\alpha r_0 \cdot  D_x^\alpha s_{-1}\big)
   +\sum_{|\alpha|=2} \frac{1}{\alpha!} \partial^\alpha_\xi r_1 \cdot D_x^\alpha s_{-1}\bigg] \nonumber \\
  && \quad \quad = \frac{1}{2}(r_1-\tau I_m)^{-1}r_1^{-1} J_m s_{-1} +s''_{-3}, \nonumber  \end{eqnarray*}
      \begin{eqnarray*} && s_{-4} =-(r_1-\tau I_m)^{-1} \bigg\{ r_0\big( \frac{1}{2} (r_1-\tau I_m)^{-1} r_1^{-1} J_m s_{-1} +s''_{-3}\big) \\    && \quad \quad \quad + \big(-\frac{1}{2} r_1^{-1} J_m +r''_{-1}\big) s_{-2}  +r_{-2} s_{-1}     \\
     && \quad \quad \quad  + \sum_{|\alpha|=1} \bigg[  \partial^\alpha_\xi r_1 \cdot D_x^\alpha\big(\frac{1}{2} (r_1-\tau I_m)^{-1}r_1^{-1}  J_m s_{-1} +s''_{-3} \big)\\
     && \quad \quad \quad +\partial^\alpha_\xi r_0 \cdot D_x^\alpha s_{-2}  +\partial^\alpha_\xi \big(-\frac{1}{2} r^{-1}_1 J_m +r''_{-1} \big)\cdot D_x^\alpha s_{-1} \bigg] \\
      && \quad \quad  \quad + \sum_{|\alpha|=2} \frac{1}{\alpha!} \bigg[\partial^\alpha_\xi r_1 \cdot D_x^\alpha s_{-2}  +\partial^\alpha_\xi r_0 \cdot D_x^\alpha s_{-1} \bigg]+
     \sum_{|\alpha|=3} \frac{1}{\alpha!} \partial^\alpha_\xi r_1\cdot D_x^\alpha s_{-1}\bigg\}\\
       && \quad \quad = -\frac{1}{2} (r_1-\tau I_m)^{-1} r_0 (r_1-\tau I_m)^{-1} r_1^{-1}  J_m s_{-1}
       +\frac{1}{2} (r_1-\tau I_m)^{-1} r_1^{-1} J_m s_{-2} \\
       && \quad \quad \quad +\frac{1}{2} (r_1-\tau I_m)^{-1} r_1^{-2} J_m r_0 s_{-1}
       +\frac{1}{4} (r_1 -\tau I_m)^{-1} r_1^{-1} \big( \partial_\xi r_1\cdot r_1^{-1} D_x J_m\big)s_{-1} \\
       && \quad \quad \quad  -\frac{E}{4}
       (r_1-\tau I_m)^{-1} r_1^{-2} J_m s_{-1} - \frac{1}{8} (r_1-\tau I_m)^{-1} r_1^{-4} \frac{\partial r_1^2}{\partial x_{n+1}} J_m s_{-1}
         \\        && \quad \quad \quad
        +\frac{1}{4} (r_1-\tau I_m)^{-1} r_1^{-2} (\partial_{x_{n+1}} J_m) s_{-1}
        \\        && \quad \quad \quad
         -\frac{1}{2} (r_1- \tau I_m)^{-1} \big(\partial_\xi r_1 \cdot (r_1 -\tau I_m)^{-1} r_1^{-1} D_x J_m\big)s_{-1}
         +s''_{-4}.\end{eqnarray*}
 Consequently, by applying the properties (\ref{-4.1})---(\ref{-4.3}) on the boundary normal coordinates once more, we have
    \begin{eqnarray*}&& r_1(x_0, \xi) =|\xi|I_m=r''_1 (x_0,\xi), \quad \;\; r_0 (x_0, \xi) =\left( \frac{1}{2|\xi|^2} \sum_{j=1}^n \kappa_j \xi_j^2 -\frac{1}{2} \sum_{j=1}^n \kappa_j\right)I_m=r''_0 (x_0, \xi),\\
 && r_{-1}(x_0, \xi)= -\frac{1}{2|\xi|}  J_m +r''_{-1}(x_0,\xi),\\
  && r_{-2} (x_0,\xi) = -\frac{1}{2|\xi|^2} \bigg( \frac{1}{2|\xi|^2} \sum_{j=1}^n  \kappa_j \xi_j^2 -\frac{1}{2} \sum_{j=1}^n \kappa_j \bigg) J_m
   -\frac{1}{4|\xi|^3}\sum_{l=1}^n \xi_l(D_{x_l} J_m)    \\
   &&\quad \quad \quad \qquad \;\; + \frac{1}{4|\xi|^2} \big(\sum_{j=1}^n \kappa_j\big)J_m +\frac{1}{4|\xi|^4} \big(\sum_{j=1}^n \kappa_j \xi_j^2 \big)J_m -\frac{1}{4|\xi|^2} \partial_{x_{n+1}}J_m +r''_{-2}(x_0,\xi)\\
  && \quad \quad \quad \quad \;\; =  \frac{1}{2|\xi|^2} \big(\sum_{j=1}^n \kappa_j \big)J_m  -\frac{1}{4|\xi|^3} \sum_{l=1}^n \xi_l (D_{x_l} J_m)
     -\frac{1}{4|\xi|^2} \partial_{x_{n+1}}J_m   +r''_{-2}(x_0,\xi). \end{eqnarray*}
so that
 \begin{eqnarray}  && s_{-1}= s''_{-1}=(|\xi|-\tau)^{-1}I_m,\nonumber \\
  && s_{-2} = s''_{-2}  =-(|\xi|-\tau)^{-2} \bigg(\frac{1}{2|\xi|^2} \sum_{j=1}^n
 \kappa_j\xi_j^2 -\frac{1}{2} \sum_{j=1}^n \kappa_j\bigg) I_m,\nonumber\\
  && s_{-3} =\frac{1}{2|\xi|} (|\xi|- \tau)^{-2}J_m +s''_{-3}, \nonumber\end{eqnarray}
  \begin{eqnarray}
   && s_{-4} = -\frac{1}{2|\xi|}  (|\xi|-\tau)^{-3} \bigg(\frac{1}{2|\xi|^2} \sum_{j=1}^n
 \kappa_j\xi_j^2 -\frac{1}{2} \sum_{j=1}^n \kappa_j\bigg) J_m \nonumber \\
 && \quad \quad \quad - \frac{1}{2|\xi|} (|\xi|-\tau)^{-3}  \bigg(\frac{1}{2|\xi|^2} \sum_{j=1}^n
 \kappa_j\xi_j^2 -\frac{1}{2} \sum_{j=1}^n \kappa_j\bigg) J_m \nonumber\\
 && \quad \quad \quad - \frac{1}{2|\xi|^2} (|\xi|-\tau)^{-2} \big( \sum_{j=1}^n
 \kappa_j\big) J_m +\frac{1}{4|\xi|^3} (|\xi|-\tau)^{-2} \sum_{l=1}^n \xi_l (D_{x_l} J_m) \nonumber \\
  && \quad \quad \quad +\frac{1}{4|\xi|^2} (|\xi|-\tau)^{-2} (\partial_{x_{n+1}}  J_m)    -\frac{1}{2|\xi|^2} (|\xi|-\tau)^{-3} \sum_{l=1}^n  \xi_l (D_{x_l}J_m) +s''_{-4}\nonumber\\
 && \quad \quad = \bigg[-\frac{1}{|\xi|} (|\xi|-\tau)^{-3}\big(\frac{1}{2|\xi|^2} \sum_{j=1}^n
 \kappa_j\xi_j^2 -\frac{1}{2} \sum_{j=1}^n \kappa_j\big) \nonumber\\
 && \quad \quad \quad -\frac{1}{2|\xi|^2} (\xi|-\tau)^{-2} (\sum_{j=1}^n \kappa_j \big)
\bigg]J_m +\frac{1}{4|\xi|^2} (|\xi|-\tau)^{-2} (\partial_{x_{n+1}}  J_m) \nonumber
\\ && \quad \quad \quad
  +\bigg( \frac{1}{4|\xi|^3} (|\xi|-\tau)^{-2} -\frac{1}{2|\xi|^2} (|\xi|-\tau)^{-3} \bigg)  \sum_{l=1}^n  \xi_l (D_{x_l}J_m)
+s''_{-4}.\nonumber\end{eqnarray}
Using the standard spherical coordinates transform, it is easy to verify (also see \cite{Liu3}) that \begin{eqnarray}  \label{-3..15}
   &&\int_{{\Bbb R}^n}  \big(\sum_{j=1}^n \xi_j^2\big)^{\frac{\beta}{2}} \xi_k^\alpha \, e^{-\sqrt{\sum_{j=1}^n \xi_j^2}}
     d\xi =\left\{ \begin{array}{ll} \Gamma (n+\beta)\,vol({\Bbb S}^{n-1})\,  &\mbox{for}\;\; \alpha=0\\
    0 \, &\mbox{for} \;\; \alpha=1,\end{array} \right. \;\; n\ge 1,\\
&&\int_{{\Bbb R}^n}   \big(\sum_{j=1}^n \xi_j^2\big)^{\frac{\beta-2}{2}} \xi_k\xi_l \, e^{-\sqrt{\sum_{j=1}^n \xi_j^2}}   d\xi =\left\{ \begin{array}{ll} \frac{\Gamma (n+\beta)\,vol({\Bbb S}^{n-1})}{n}   \, &\mbox{for}\;\; k=l\\
0 \;\; &\mbox{for}\;\; k\ne l, \end{array} \right.\;\; n\ge 2, \nonumber\end{eqnarray}
 The following Cauchy's differentiation formula will be used, which states that for any complex analytic function $f$.
\begin{eqnarray} \label{-3..16} \frac{i}{2\pi} \int_\Gamma \frac{f(\tau)}{(z-\tau)^{k+1}} d\tau =\frac{(-1)^k}{k!} f^{(k)} (z).\end{eqnarray}
In view of (\ref{3./1}), (\ref{-3..15}) and (\ref{-3..16}), we find that
\begin{eqnarray} && \label{3---14}a_{m,l}(x,n)= \frac{1}{(2\pi)^n}
 \int_{{\Bbb R}^n} \bigg( \frac{i}{2\pi}  \int_\Gamma  e^{-\tau }  s_{-1-l} (x, \xi, \tau) d\tau\bigg) d\xi \qquad \qquad \\
&& \quad \quad \quad \, =\frac{1}{(2\pi)^n}
 \int_{{\Bbb R}^n} \bigg( \frac{i}{2\pi}  \int_\Gamma e^{-\tau }  s''_{-1-l} (x, \xi, \tau) d\tau\bigg) d\xi =a''_{1,l}(x,n),\, \quad l=0,1,
  \nonumber\end{eqnarray}
    \begin{eqnarray}\label{3--16} && a_{m,2}(x,n)= \frac{1}{(2\pi)^n}
 \int_{{\Bbb R}^n} \bigg( \frac{i}{2\pi}  \int_\Gamma  e^{-\tau }  s_{-3} (x, \xi, \tau) d\tau\bigg) d\xi  \\
  && \quad \quad\quad\quad \;\;= \frac{1}{(2\pi)^n}
 \int_{{\Bbb R}^n} \bigg( \frac{i}{2\pi}  \int_\Gamma  e^{-\tau } \big(\frac{1}{2|\xi|} (|\xi|-\tau)^{-2} J_m +s''_{-3} \big) d\tau \bigg)
 d\xi \nonumber\\
  &&  \quad \quad\quad\quad \;\;= \frac{1}{(2\pi)^n}
 \int_{{\Bbb R}^n} \frac{1}{2|\xi|} e^{-|\xi| } J_m \,d\xi +\frac{1}{(2\pi)^n} \int_{{\Bbb R}^n} \bigg(\frac{i}{2\pi} \int_\Gamma e^{-\tau }
   s''_{-3} d\tau  \bigg) d\xi \nonumber\\
    &&  \quad \quad\quad\quad \;\;=  \frac{ \Gamma(n-1) \cdot vol({\Bbb S}^{n-1})}{2(2\pi)^n}\, J_m  +a''_{1,2}(x,n),\nonumber\end{eqnarray}
\begin{eqnarray} \label{3---17} && \quad \quad a_{m,3}(x,n)= \frac{1}{(2\pi)^n}
 \int_{{\Bbb R}^n} \bigg( \frac{i}{2\pi}  \int_\Gamma  e^{-\tau }  s_{-4} (x, \xi, \tau) d\tau\bigg) d\xi  \\
  && \quad \quad\quad\quad \;\, =\frac{1}{(2\pi)^n}
 \int_{{\Bbb R}^n}  \frac{i}{2\pi}  \int_\Gamma  e^{-\tau }\bigg[ \bigg(\frac{1}{2|\xi|} \sum_{j=1}^n \kappa_j
      - \frac{1}{2|\xi|^3}  \sum_{j=1}^n \kappa_j \xi_j^2 \bigg) (|\xi|-\tau)^{-3} J_m \nonumber\\
    && \quad \quad\quad\quad \quad\; -  \frac{1}{2|\xi|^2} (|\xi|-\tau)^{-2}\big(\sum_{j=1}^n \kappa_j \big) J_m
    +\frac{1}{4|\xi|^2} (|\xi|-\tau)^{-2} (\partial_{x_{n+1}} J_m) \nonumber \\ && \quad \quad\quad\quad \quad\; +\bigg( \frac{1}{4|\xi|^3} (|\xi|-\tau)^{-2}
    -\frac{1}{2|\xi|^2} (|\xi|-\tau)^{-3}\bigg) \sum_{l=1}^n  \xi_l ( D_{x_l} J_m)
        + s''_{-4} \bigg] d\tau\,d\xi\nonumber\\
      && \quad \quad\quad\quad \;\, = \frac{1}{(2\pi)^n}
 \int_{{\Bbb R}^n} e^{-|\xi|}  \bigg[\bigg(\frac{1}{4|\xi|} \sum_{j=1}^n \kappa_j
      - \frac{1}{4|\xi|^3}  \sum_{j=1}^n \kappa_j \xi_j^2  -\frac{1}{2|\xi|^2} \sum_{j=1}^n \kappa_j  \bigg) J_m \nonumber\\
      && \quad \quad\quad\quad \quad\;
      +\frac{1}{4|\xi|^2} (\partial_{x_{n+1}}J_m) +\big( \frac{1}{4|\xi|^3}
    -\frac{1}{4|\xi|^2} \big) \sum_{l=1}^n  \xi_l (D_{x_l} J_m)
      \bigg]\, d\xi \nonumber\\
      && \quad \quad\quad\quad \quad\;
      +\frac{1}{(2\pi)^n} \int_{{\Bbb R}^n} \bigg( \frac{i}{2\pi} \int_\Gamma e^{-\tau } s''_{-4} d\tau\bigg) d\xi\nonumber\\
            && \quad \quad\quad\quad \;\,= \frac{1}{(2\pi)^n}\bigg[ \bigg(\frac{n-1}{4n} \Gamma(n-1)   -\frac{1}{2} \Gamma(n-2) \bigg) \big(\sum_{j=1}^n  \kappa_j\big)J_m \nonumber\\
      && \quad \quad\quad\quad \quad\;      + \frac{\Gamma(n-2)}{4}  (\partial_{x_{n+1}} J_m)\bigg]vol({\Bbb S}^{n-1})
            +a''_{1,3}(x,n).\nonumber \end{eqnarray}
By comparing the symbol equations (1.5) of \cite{LU} (or, the symbol expression on p.$\,$1103 of \cite{LU}) with (\ref{2-7}), and by applying theorem 6.1 of \cite{Liu3} we immediately see that
  \begin{eqnarray} \label{06--060-1} a''_{1,0}(n,x)=\frac{\Gamma(\frac{n+1}{2})}{\pi^{\frac{n+1}{2}}}I_m,\end{eqnarray}
  \begin{eqnarray}\label{06--060}\quad a''_{1,1}(n,x)= \left(\frac{1}{2\pi}\right)^{n}   \frac{(n-1)\Gamma(n)\,\mbox{vol}({\Bbb S}^{n-1})}{2n} \bigg(\sum_{j=1}^n \kappa_j(x)\bigg)I_m,
     \end{eqnarray}
    \begin{eqnarray}  \label{66..55}&&    a''_{1,2}(n,x) = \frac{\Gamma(n-1) \mbox{vol}({\Bbb S}^{n-1}) }{8(2\pi)^n} \left[ \frac{3-n}{3n} R_{\partial \Omega} + \frac{n-1}{n} {\tilde R}_{\Omega}\right. \\ && \left. +\frac{n^3 -n^2 -4n +6}{n(n+2)} \big(\sum_{j=1}^n \kappa_j (x)\big)^2 + \frac{n^2 -n-2}{n(n+2)} \sum_{j=1}^n \kappa_j^2(x) \right] I_m
      \nonumber \end{eqnarray} and
\begin{eqnarray} \label{6.6--1.} && a''_{1,3}(n,x)=
  \bigg(\frac{1}{2\pi}\bigg)^{n}\frac{\Gamma(n-2)\,\mbox{vol} (S^{n-1}) }{8n}  \bigg[
   \frac{n^3 -2n^2 -7n+7}{2(n+2)} {\tilde R}_\Omega(x)\, \big(\sum_{j=1}^n \kappa_j(x)\big) \\ && \quad \, + \frac{-3n^4 -4 n^3 +59n^2 +75 n -180}{ 6(n+2)(n+4)} \big(\sum_{j=1}^n \kappa_j (x) \big)R_{\partial \Omega} \nonumber\\ && \quad \,+
   \frac{n^5 -20 n^3 +2n^2 +61 n -74}{6(n+2)(n+4)} \big(\sum_{j=1}^n \kappa_j(x)\big)^3 \nonumber\\ && \quad \,
    + \frac{n^4 +8n^3 +15n^2 +3n -32}{2(n+2) (n+4)} \big(\sum_{j=1}^n \kappa_j(x)\big)\big(\sum_{j=1}^n \kappa_j^2 (x)\big)\nonumber\\ && \quad \,
    +\frac{-6n^3-34n^2 +40}{3(n+2)(n+4)} \sum_{j=1}^n \kappa_j^3(x) + \frac{4n^2-6}{n+2} \sum_{j=1}^n\kappa_j(x) {\tilde R}_{jj}(x)\nonumber
    \\ && \quad \,-\frac{12n^3 +50 n^2 -6n -104}{3(n+2)(n+4)} \sum_{j=1}^n \kappa_j(x) R_{jj}(x) + (n-1) \sum_{j=1}^n {\tilde{R}}_{j(n+1)j(n+1),(n+1)}(x)\nonumber\\ &&\quad \,-\frac{n-2}{2 } {\tilde{R}}_{\Omega} (x)  +\frac{n-2}{2} R_{\partial \Omega} -\frac{n-2}{2} \sum_{j=1}^n \kappa_j^2(x)\bigg]I_m \nonumber
.\end{eqnarray}
Inserting (\ref{06--060-1})---(\ref{6.6--1.}) into (\ref{3---14})---(\ref{3---17}) we get the desired result for $a_{m,l} (x,n)$, $\; (m\ge 1,\;\, 0\le l\le 3)$. $\quad \ \square$

\vskip 0.28 true cm

\noindent{\bf Remark 3.2.} \ \  (i) \  Our asymptotic expansion is sharp because the  asymptotic expansion (\ref{1,,1}) holds  for $t\to 0^+$ (see, (4.2.62) of \cite{Gr}). In particular, when $m=1$ and $q\equiv 0$ we recover the results in \cite{Liu3} for $a_0, a_1, a_3, a_4$, and in \cite{PS} for $a_0,a_1,a_2$.

 (ii) \   By our method, we can also get the asymptotic expansion for any integer $M>4$:
\begin{eqnarray} &&   \int_{\partial \Omega}  {\mathcal{K}}_m(t,x,x) dS(x)= t^{-n}\frac{\Gamma(\frac{n+1}{2})}{\pi^{\frac{n+1}{2}}} \int_{\partial \Omega} I_m\, dS(x) + t^{1-n} \int_{\partial \Omega} a_{m,1}(n,x)\,dS(x)\\
  && \qquad\quad \quad +t^{2-n} \int_{\partial \Omega} a_{m,2}(n,x)\, dS(x)   +\cdots+ t^{M-1-n}  \int_{\partial \Omega}
 a_{m,M-1}(n,x)\, dS(x) \nonumber\\
  && \qquad\quad \quad+ \left\{\begin{array}{ll} O(t^{M-n})I_m \;\quad \, \mbox{when}\;\; n>M-1,\nonumber\\ O(t\log t)I_m \;\quad\, \mbox{when} \;\; n=M-1,\end{array}\right.\quad \;\;\, \mbox{as}\;\; t\to 0^+, \nonumber \end{eqnarray}
if we further calculate the lower-order symbol equations for the operator $\Lambda_m$.

(iii) \ \ The above \ asymptotic
expansion  \ shows that one can ``hear'' \ all heat invariants  $\int_{\partial \Omega} a_{M-1}(n,x)dS(x)$ ($M=1,2,3,\cdots$)  by ``hearing'' all of the perturbed polyharmonic Steklov eigenvalues.

(iv) \  By applying the Tauberian theorem (see, for example, Theorem 15.3 of p.$\,$30 of \cite{Kor}) for the first term on the right side of (\ref{6.0.1}), we immediately get the Weyl-type law  for the perturbed polyharmonic Steklov eigenvalues (when $m=1$ and $q\equiv 0$, it is just Sandgren's asymptotic formula, see \cite{Sa}):
\begin{eqnarray} \label{-4.13}  N(\lambda)=
 \#\{k\big|\lambda_k \le \lambda\} =\left[
   \frac{\omega_{n}\big(\mbox{vol}(\partial \Omega)\big)\lambda^{n}}{(2\pi)^{n}}
  +o(\lambda^{n})\right] \quad \;\mbox{as}\;\;
 \lambda\to +\infty.\end{eqnarray}

\vskip 0.98 true cm

\centerline {\bf  Acknowledgments}

\vskip 0.38 true cm
 I wish to express my sincere gratitude to Professor
 L. Nirenberg and Professor Fang-Hua Lin for their
  support and help.  This
research was supported by SRF for ROCS, SEM (No. 2004307D01)
   and NNSF of China (11171023/A010801).

  \vskip 1.65 true cm

\vskip 0.32 true cm

\end{document}